\documentclass[a4paper,11pt,reqno]{amsart}
\usepackage{amssymb,amsthm,amsmath}
\usepackage{mathtools}
\usepackage{amsfonts}
\usepackage{amscd}
\usepackage{amssymb}
\usepackage{enumerate}
\usepackage{bbm}
\usepackage{graphicx}
\usepackage{mathabx}
\usepackage{color}
\usepackage{amsbsy}
\usepackage{graphicx}
\usepackage{amsthm}
\usepackage{amsmath}
\usepackage{amsxtra}
\usepackage{mathrsfs}
\usepackage{bbm}
\usepackage{dsfont}
\usepackage{enumitem}
\usepackage{comment}
\usepackage{enumerate}
\usepackage{xcolor}
\usepackage{float}
\usepackage{url}
\usepackage{soul}
\usepackage[hidelinks]{hyperref} 

\allowdisplaybreaks

\usepackage{ifthen}

\hypersetup{
    colorlinks=true,
    linkcolor=red,      
    citecolor=red,       
    urlcolor=cyan        
}

\newtheorem{theorem}{Theorem}[section]
\newtheorem{thm}{Theorem}[subsection]

\newtheorem{cor}[thm]{Corollary}

\newtheorem{lem}[thm]{Lemma}

\numberwithin{equation}{subsection}
\definecolor{red}{rgb}{1.0, 0.0, 0.0}
\newcommand{\Bea}{\begin{eqnarray*}}
	\newcommand{\Eea}{\end{eqnarray*}}
\newcommand{\Be} {\begin{equation*}}
	\newcommand{\Ee} {\end{equation*}}
\newcommand{\be} {\begin{equation}}
	\newcommand{\ee} {\end{equation}}
\newcommand{\bea} {\begin{eqnarray}}
	\newcommand{\eea} {\end{eqnarray}}

\usepackage{color}

\title[$L^p$-HPW Uncertainty Inequalities on the Laguerre Hypergroup]{$L^p$-Heisenberg--Pauli--Weyl Uncertainty Inequalities on the Laguerre Hypergroup}

\author[A.~Dabra]{Arvish Dabra}
\address{
    Arvish Dabra:
    \endgraf
    Department of Mathematics
    \endgraf
    Indian Institute of Technology Delhi 
    \endgraf
    Hauz Khas, New Delhi - 110016
    \endgraf
    India
    \endgraf
    {\it E-mail address:} {\rm arvishdabra3@gmail.com}
}

\author[A.~Dasgupta]{Aparajita Dasgupta}
\address{
    Aparajita Dasgupta:
    \endgraf
    Department of Mathematics
    \endgraf
    Indian Institute of Technology Delhi 
    \endgraf
    Hauz Khas, New Delhi - 110016
    \endgraf
    India
    \endgraf
    {\it E-mail address:} {\rm adasgupta@maths.iitd.ac.in}
}


    \vspace{1cm}

    \allowdisplaybreaks
\begin{document}
    
\begin{abstract}
In this paper, we establish the first $L^p$-Heisenberg--Pauli--Weyl uncertainty inequalities on the Laguerre hypergroup for the full range $1\le p\le2$. These results extend Xiao's Euclidean $L^p$ theory to the setting of the Laguerre hypergroup, which is the fundamental manifold of the radial function space for the Heisenberg group. The analysis is carried out through the Fourier--Laguerre transform and exploits the mixed discrete--continuous spectral structure of the Laguerre hypergroup, requiring estimates adapted to its Plancherel measure and dilation structure. As a consequence, in the endpoint case $p=2$, we obtain a refined $L^2$-Heisenberg--Pauli--Weyl uncertainty inequality valid for all positive exponents $a,b>0$, thereby improving the earlier result of Atef (2013), where the assumptions $a,b\ge1$ arose from the heat kernel methods. Our proofs rely on the Fourier--Laguerre transform, dilation and scaling invariance, the Hausdorff--Young inequality and the Plancherel identity, completely avoiding heat kernel techniques. These results provide a unified Fourier-analytic framework for Heisenberg--Pauli--Weyl uncertainty inequalities on the Laguerre hypergroup and further strengthen the connections between Euclidean, Heisenberg and hypergroup harmonic analysis.
\end{abstract}
   
  \maketitle
  

    \section{Introduction}

    The uncertainty principle is a fundamental concept in harmonic analysis, asserting that \emph{a non-zero function and its Fourier transform cannot both be sharply localized}. In simpler terms, if a function is \emph{highly concentrated}, its Fourier transform cannot also be, unless the function is identically zero. This concept was first articulated in a mathematical framework by Heisenberg in his seminal 1927 paper \cite{Heisen} and was later developed more formally by Weyl \cite{weyl}. The classical Heisenberg--Pauli--Weyl (HPW) uncertainty inequality \cite[(3.1)]{FSsurvey} for the Fourier transform on $\mathbb{R}^n$ is
        $$ \|f\|_2^4 \leq C_n \, \left( \int_{\mathbb{R}^n} |x|^2 \, |f(x)|^2 \, dx \right) \, \left( \int_{\mathbb{R}^n} |\xi|^2 \, |\hat{f}(\xi)|^2 \, d\xi \right).$$
    By the Plancherel identity, this inequality is equivalently expressed in terms of the Laplacian $\Delta$ on $\mathbb{R}^n$ as
    $$\|f\|_2^2 \leq C_n \, \| \, |x| \, f\|_2 \, \|(-\Delta)^{1/2} \, f\|_2.$$
    Since its inception, the uncertainty principle has been extensively studied in harmonic analysis, mathematical physics, partial differential equations and signal processing. For a comprehensive overview of its history, significance and numerous variants, we refer the reader to the survey article of Folland and Sitaram \cite{FSsurvey}.

    Motivated by the Euclidean theory, uncertainty inequalities have been extensively investigated in noncommutative settings, particularly on the Heisenberg group $\mathbb{H}^n$. In 1990, Thangavelu \cite{Th1990} established the following HPW inequality for normalized functions in $L^2(\mathbb{H}^n)$:
    $$\left( \int_\mathbb{R} \int_{\mathbb{C}^n} |(z-a)|^2 \, |f(z,t)|^2 \, dz \, dt \right)^{1/2} \, \left( \int_\mathbb{R} |\lambda|^n \, \|\hat{f}(\lambda) (H(\lambda))^{1/2} \|_{HS}^2 \, d\lambda \right)^{1/2} \geq \sqrt{n} \, \left(\frac{\pi}{2}\right)^{\frac{n+1}{2}}.$$
    Subsequently, Sitaram et al.~\cite{SST1995} generalized this result by establishing the following uncertainty inequality involving the Heisenberg sub-Laplacian $\mathcal{L}.$
    \begin{theorem}\label{1.1}
        For $f \in L^2(\mathbb{H}^n)$ and $0 \leq \gamma < Q/2 = n+1,$ we have
        $$\|f\|_2^2 \leq C \, \left( \int_{\mathbb{H}^n} |(z,t)|^{2\gamma} \, |f(z,t)|^2 \, dz \, dt \right)^{1/2} \, \left( \int_{\mathbb{H}^n} |\mathcal{L}^{\gamma/2} \, f(z,t)|^2 \, dz \, dt \right)^{1/2},$$
    where $C$ is a positive constant independent of $f$.
    \end{theorem}
    In 2012, Xiao and He \cite{XiaoHe} further extended Theorem \ref{1.1} by removing the restriction $0 \leq \gamma < Q/2$, thereby establishing a more general HPW uncertainty inequality on the Heisenberg group.\\
  Over the past two decades, HPW-type uncertainty inequalities have been established in a variety of non-Euclidean settings. Ciatti et al.~\cite{CRSKMA} proved HPW uncertainty inequalities on two-step nilpotent Lie groups and subsequently extended these results to Lie groups with polynomial volume growth \cite{CRS2007}. Ruzhansky and Suragan \cite[9.3.1]{RuzDur} established an HPW-type uncertainty principle on homogeneous groups of homogeneous dimension $Q \geq 3$. Furthermore, Ciatti et al.~\cite{CCR20215} obtained the following $L^p$ version of the HPW uncertainty inequality on a stratified Lie group $\mathbb{G}$ of homogeneous dimension $Q$ endowed with a positive hypoelliptic sub-Laplacian $\mathcal{L}$.
    \begin{theorem}
        Suppose that $|\cdot|$ is a homogeneous norm on $\mathbb{G}$ such that $\beta,\delta >0,$ $p,r >1$ and $s \geq 1$ with 
        $$\frac{\beta+\delta}{p} = \frac{\delta}{s} + \frac{\beta}{r}.$$
        Then
        $$\|f\|_p \leq C \, \|\,|\cdot|^\beta \, f\|_s^{\frac{\delta}{\beta + \delta}} \, \|\mathcal{L}^{\delta/2} \, f\|_r^{\frac{\beta}{\beta+\delta}}$$
        for every $f \in \mathcal{C}_c^\infty(\mathbb{G}).$
    \end{theorem}
    In the Euclidean setting, Steinerberger \cite{SStei} established an $L^1$ version of the HPW uncertainty inequality on $\mathbb{R}^n$. More recently, Xiao \cite{XiaoRn} developed a general $L^p$ theory by proving the $L^p$-Heisenberg--Pauli--Weyl uncertainty inequalities for the full range $1 \leq p < \infty$. More precisely, the following result holds.
    \begin{theorem}\label{1.3}
        Let $1 \leq p < \infty$ and $1 < p' = p/(p-1) \leq \infty.$ Let $0 < a,b < \infty$ such that $b > n(1/p - 1/2)$ whenever $1 \leq p < 2$ and $a < n/p$ whenever $p \geq 2.$ Then, for every $f \in L^p(\mathbb{R}^n),$ we have
        $$\|f\|_p^{a+b} \leq C(a,b,p,n) \, \| \, |x|^a \, f\|_p^b \, \| \, |\xi|^b \, \hat{f}\|_{p'}^a,$$
        where $C(a,b,p,n)$ is a positive constant depending only on the parameters.
    \end{theorem}
     In 2023, Fu and Xiao \cite{FX} extended this line of research by proving an uncertainty principle within the framework of Lorentz spaces. In 2025, Ganguly and Sarkar \cite{GS} proved an analogue of Theorem \ref{1.3} on two-step MW groups for $1 \leq p \leq 2$, while Chen et al.~\cite{CDM} established $L^p$-type HPW uncertainty principles for the fractional Fourier transform.\\
 Motivated by the works of Steinerberger \cite{SStei}, Xiao \cite{XiaoRn} and Ganguly and Sarkar \cite{GS}, we investigate the $L^p$-Heisenberg--Pauli--Weyl uncertainty inequalities on the Laguerre hypergroup, which is the fundamental manifold of the radial function space for the Heisenberg group. Since the Laguerre hypergroup encodes the radial harmonic analysis of the Heisenberg group, the uncertainty inequalities established in this paper may be viewed as radial analogues of the corresponding inequalities on the Heisenberg group. Consequently, the present results provide insight into the interaction between commutative hypergroup analysis and noncommutative harmonic analysis.\\
The Laguerre hypergroup is a fundamental example of a commutative hypergroup and provides a  natural structure that allows one to extend Fourier and spectral analysis techniques to non-Euclidean settings. More generally, the concept of a \emph{hypergroup} was introduced by Jewett (1975) \cite{Je} as a generalization of locally compact groups by replacing the product of two elements with a probability measure. This theory has found numerous applications in harmonic analysis, probability theory and special functions. Recent developments include the work of Kumar and Ruzhansky \cite{VR2022} on Hardy--Littlewood inequalities and $L^p-L^q$ Fourier multipliers on compact and commutative hypergroups \cite{VM}. 
    
    The Laguerre hypergroup is defined as
\[
\mathbb{K} = [0,\infty) \times \mathbb{R},
\]
equipped with the Radon measure
\[
dm_{\alpha}(x,t) = \frac{1}{\pi \, \Gamma(\alpha+1)} \, x^{2\alpha+1} dx \, dt,
\]
and  the convolution structure based on Laguerre functions and spherical averaging. The corresponding generalized sub-Laplacian
\[
\mathcal{L} = -\left( \frac{\partial^2}{\partial x^2} + \frac{2\alpha+1}{x}\frac{\partial}{\partial x} + x^2 \frac{\partial^2}{\partial t^2} \right),
\]
is a positive self-adjoint operator on $L^2(\mathbb{K})$. When $\alpha=n-1$, it coincides with the radial part of the sub-Laplacian on the Heisenberg group $\mathbb{H}^n$ \cite{KS}. The associated Fourier--Laguerre transform provides an isometric isomorphism between $L^2(\mathbb{K})$ and $L^2(\widehat{\mathbb{K}})$, where $\widehat{\mathbb{K}} = \mathbb{R}\times \mathbb{N}$. This yields hypergroup analogues of the classical Fourier analysis results such as the Plancherel identity, the Hausdorff--Young inequality and inversion formula \cite{NT}. Further details on the Laguerre hypergroup is provided in Section \ref{sec2}. \\
Over the past two decades, harmonic analysis on the Laguerre hypergroup has been the subject of extensive research. Key contributions include pseudo-differential operators \cite{MAssal}, Sobolev-type spaces on the dual \cite{AN2004}, maximal functions \cite{gulass} and fractional integral inequalities \cite{EH2017,GO}. In the context of uncertainty principles, Huang and Liu \cite{HuLiu,HuLiu2011} established analogues of Beurling’s and Hardy’s theorems. Subsequently, Atef (2013) \cite{AR} proved the following $L^2$-HPW uncertainty inequality on the Laguerre hypergroup using heat kernel techniques.
    
    \begin{theorem}\label{1.4}
        Let $a,b \geq 1.$ Then, for every $f \in L^2(\mathbb{K}),$ we have
        $$\| f \|_2^{a+b} \lesssim \| \, |(x,t)|^a \, f\|_2^b \, \| \, |(\lambda,m)|^{b/2} \, \hat{f}\|_2^a.$$
    \end{theorem}
 
    More recently, Tyr and Daher \cite{Tyr2024, TD2022} investigated Bernstein's and Jackson's inequalities in this framework.

    Thus, the Laguerre hypergroup serves as a bridge between Euclidean analysis, Lie group harmonic analysis and abstract hypergroup theory. Motivated by the above developments, we establish new $L^p$-Heisenberg--Pauli--Weyl uncertainty inequalities on the Laguerre hypergroup $\mathbb{K}$ of homogeneous dimension $Q = 2\alpha+4$. Our results extend Xiao's Euclidean $L^p$-HPW uncertainty inequalities to the Laguerre hypergroup setting and refine the $L^2$-HPW uncertainty inequality established by Atef \cite{AR}.
    
    More precisely, we establish the following $L^p$-HPW uncertainty inequality on the Laguerre hypergroup for the full range $1 \leq p \leq 2$. 
    
    \begin{theorem}\label{1.5}
        Let $1 \leq p \leq 2$ and $2 \leq p' = p/(p-1) \leq \infty.$ Suppose that $a > 0$ and $b > Q(1/p - 1/2).$ Then, for every $f \in L^p(\mathbb{K}),$ we have
            $$\|f\|_p^{a+b} \leq C_{a,b,p,\alpha} \, \| \, |(x,t)|^a \, f\|_p^b \, \| \, |(\lambda,m)|^{b/2} \, \hat{f}\|_{p'}^a,$$
        where $C_{a,b,p,\alpha}$ is a positive constant depending only on the parameters.
    \end{theorem}

    In contrast to the Euclidean and stratified Lie group settings, harmonic analysis on the Laguerre hypergroup is carried out through the Fourier--Laguerre transform, whose spectral side is parameterized by the mixed discrete--continuous variables $(\lambda,m)$, equipped with the Plancherel measure $d\gamma_\alpha(\lambda,m)$. Consequently, several key components of the proof require estimates adapted to this spectral structure, including the analysis of the Plancherel measure of spectral balls and the treatment of weighted Fourier norms. Thus, although our strategy is inspired by the Euclidean work of Xiao and the analysis of MW groups due to Ganguly and Sarkar, its implementation requires nontrivial modifications dictated by the hypergroup setting.
    
    As a notable consequence, for $p = 2$, we obtain the following \emph{refined} classical Heisenberg--Pauli--Weyl uncertainty inequality.

    \begin{theorem}\label{1.7}
        Let $a,b > 0$. Then, for every $f \in L^2(\mathbb{K}),$ we have
        $$\| f \|_2^{a+b} \lesssim \| \, |(x,t)|^a \, f\|_2^b \, \| \, |(\lambda,m)|^{b/2} \, \hat{f}\|_2^a.$$
    \end{theorem}
    
   Theorem \ref{1.7} refines Atef's $L^2$-HPW uncertainty inequality (Theorem \ref{1.4}) by removing the restrictive assumptions $a,b \geq 1$ and establishing the inequality for all positive constants $a$ and $b$. Furthermore, the proof avoids the use of heat kernel techniques.
  
    The present work establishes a unified $L^p$ theory of Heisenberg--Pauli--Weyl uncertainty inequalities on the Laguerre hypergroup for the full range $1 \leq p \leq 2$. These results extend Xiao's Euclidean $L^p$-HPW uncertainty inequalities \cite{XiaoRn} to the Laguerre hypergroup setting and complement the recent $L^p$ theory developed by Ganguly and Sarkar \cite{GS} for two-step MW groups. Since the Laguerre hypergroup arises naturally in radial harmonic analysis on the Heisenberg group, the present results place the Euclidean, nilpotent Lie group and Laguerre hypergroup settings within a common $L^p$ framework for Heisenberg--Paul--Weyl uncertainty inequalities. This work marks a significant contribution to the study of uncertainty principles in non-Euclidean harmonic analysis.
    
    Besides establishing the first $L^p$-Heisenberg--Pauli--Weyl uncertainty inequalities on the Laguerre hypergroup for the full range $1\le p\le2$, our results also provide a genuine improvement of the earlier $L^2$-result of Atef \cite{AR}. In particular, Atef's theorem was proved only under the assumptions $a,b\ge1$, which arise from the heat kernel method employed in that work. In contrast, our approach is entirely Fourier-analytic, relying on the Fourier--Laguerre transform, dilation and scaling invariance, the Hausdorff--Young inequality and the Plancherel identity. We first establish the uncertainty inequality in the natural parameter range $0<a<Q/2$ and $0<b<4$ and then use H\"older interpolation argument to extend it to all positive exponents $a,b>0$. Consequently, the restrictions $a,b\ge1$ are shown to be artifacts of the previous proof technique rather than intrinsic limitations of the uncertainty inequality itself. This yields a sharper formulation of the classical $L^2$-Heisenberg--Pauli--Weyl uncertainty principle on the Laguerre hypergroup and provides a proof strategy that is both conceptually simpler and potentially applicable to other hypergroup settings. We emphasize, however, that our objective is not to determine the optimal constants in these inequalities. Since our proofs rely on interpolation and Fourier-analytic techniques rather than extremal variational methods, the sharp constants and the corresponding extremizers remain interesting open problems.

 The paper is organized as follows:   In Section~\ref{sec2}, we recall the necessary background on the Laguerre hypergroup. 
    We begin by presenting the definition of the Laguerre hypergroup and the associated Radon measure, followed by a description of the hypergroup convolution and involution. We also discuss the generalized sub-Laplacian, its relation to the Heisenberg group and the construction of the Fourier--Laguerre transform, together with fundamental results such as the Plancherel identity, inversion formula and scaling properties. These preliminaries provide the analytic framework necessary for proving our uncertainty inequalities.
 In Section~\ref{sec3}, we present the proof of our main result, Theorem \ref{1.5}. The analysis is divided into three cases according to the range of $p$: 

\begin{enumerate}
    \item \textbf{Case $p=1$.} We establish the inequality by exploiting dilation invariance together with sharp $L^1$ estimates. 
    \item \textbf{Case $1 < p < 2$.} The proof relies on interpolation techniques, the Plancherel identity and precise control of the Fourier--Laguerre transform. 
    \item \textbf{Case $p = 2$.} The argument employs auxiliary functions, $L^2$ bounds for exponential multipliers and a careful balance between spatial and frequency localization.
\end{enumerate}
Finally, these results are combined to yield Theorem~\ref{1.5}, which synthesizes the $L^p$-Heisenberg--Pauli--Weyl uncertainty inequalities on the Laguerre hypergroup for the full range $1 \leq p \leq 2$.
\section{Preliminaries and Notations}\label{sec2}
This section collects the analytic background and notation required for our study of uncertainty inequalities on Laguerre hypergroup $\mathbb{K}$. 
In particular, we recall the measure structure, convolution, Fourier--Laguerre transform and dilation properties of $\mathbb{K}$. 
These tools form the analytic framework within which the proof of our main result will be carried out. 
For more general references on hypergroups, see \cite{BH, Je, Tri}.

\medskip

Let $\alpha \geq 0$. We define 
\[
\mathbb{K} = [0,\infty) \times \mathbb{R},
\]
equipped with the Radon measure
\[
dm_\alpha(x,t) = \frac{1}{\pi \, \Gamma(\alpha+1)} \, x^{2\alpha+1} \, dx \, dt.
\]
Denote by $M(\mathbb{K})$ the space of bounded Radon measures on $\mathbb{K}$. For $\mu,\nu \in M(\mathbb{K})$, the convolution $\mu \ast \nu$ is defined by
\[
(\mu \ast \nu)(f) := \int_{\mathbb{K} \times \mathbb{K}} T^{(\alpha)}_{(x,t)} f(y,s) \, d\mu(x,t)\, d\nu(y,s),
\]
where the translation operator $T^{(\alpha)}_{(x,t)}$ is given by
{\small{\[
T^{(\alpha)}_{(x,t)} f(y,s) =
\begin{dcases}
\frac{1}{2\pi}\int_0^{2\pi} 
f\!\left(\sqrt{x^2+y^2+2xy\cos\theta},\, s+t+xy\sin\theta\right)\, d\theta, & \alpha=0, \\[1.0em]
\frac{\alpha}{\pi} \int_0^{2\pi} \int_0^1 
f\!\left(\sqrt{x^2+y^2+2xyr\cos\theta},\, s+t+xyr\sin\theta\right)\, r(1-r^2)^{\alpha-1}\, dr\, d\theta, & \alpha >0.
\end{dcases}
\]}}
It follows that $(\mathbb{K},\ast,i)$ is a commutative hypergroup in the sense of Jewett \cite{Je}, with involution
\[
i(x,t) = (x,-t).
\]
A natural homogeneous norm on $\mathbb{K}$ is
\[
|(x,t)| := \left(x^4+4t^2\right)^{1/4}.
\]

For $1\leq p \leq \infty$, the Lebesgue space $L^p(\mathbb{K})$ consists of measurable functions $f:\mathbb{K}\to \mathbb{C}$ with finite norm
\[
\|f\|_p :=
\begin{dcases}
\left( \int_{\mathbb{K}} |f(x,t)|^p \, dm_\alpha(x,t)\right)^{1/p}, & 1\leq p <\infty, \\[0.8em]
\operatorname*{ess\,sup}_{(x,t)\in \mathbb{K}} |f(x,t)|, & p=\infty.
\end{dcases}
\]

The Laguerre hypergroup carries a distinguished differential operator,
\[
\mathcal{L} = -\left(\frac{\partial^2}{\partial x^2} + \frac{2\alpha+1}{x}\frac{\partial}{\partial x} + x^2 \frac{\partial^2}{\partial t^2}\right),
\]
which is non-negative and symmetric on $L^2(\mathbb{K})$. When $\alpha=n-1$, this operator coincides with the radial part of the sub-Laplacian on the Heisenberg group $\mathbb{H}^n$ \cite{KS, KS90, stempak2026laguerre}.

For each $(\lambda,m)\in \mathbb{R}\times \mathbb{N}$, the Laguerre function
\[
\varphi_{(\lambda,m)}(x,t) = \frac{m!\,\Gamma(\alpha+1)}{\Gamma(m+\alpha+1)} \, e^{i\lambda t}\, e^{-|\lambda|x^2/2} \, L_m^\alpha(|\lambda|x^2),
\]
is an eigenfunction of the operator $\mathcal{L}$ associated with the eigenvalue
$$4|\lambda|\left(m+\frac{\alpha+1}{2}\right),$$
that is,
\[
\mathcal{L}\varphi_{(\lambda,m)} = 4|\lambda|\left(m+\frac{\alpha+1}{2}\right)\,\varphi_{(\lambda,m)}.
\]

For $f \in L^1(\mathbb{K})$, the Fourier--Laguerre transform is defined by
\[
\widehat{f}(\lambda,m) := \int_{\mathbb{K}} f(x,t)\, \varphi_{(-\lambda,m)}(x,t)\, dm_\alpha(x,t).
\]
The dual object $\widehat{\mathbb{K}} = \mathbb{R}\times \mathbb{N}$ is equipped with the Plancherel measure
\[
d\gamma_\alpha(\lambda,m) = L_m^\alpha(0)\, |\lambda|^{\alpha+1}\, d\lambda,
\]
so that
\[
\int_{\widehat{\mathbb{K}}} g(\lambda,m)\, d\gamma_\alpha(\lambda,m) = \sum_{m\geq 0} L_m^\alpha(0) \int_{\mathbb{R}} g(\lambda,m)\, |\lambda|^{\alpha+1}\, d\lambda.
\]
The quasi-norm on $\mathbb{\widehat{K}}$ is given by
\[
|(\lambda,m)|:= 4 |\lambda|\left(m + \frac{\alpha+1}{2}\right).
\]
The Fourier--Laguerre transform extends to an isometric isomorphism $\mathcal{F}:L^2(\mathbb{K})\to L^2(\widehat{\mathbb{K}})$ and satisfies the Plancherel identity
\[
\|f\|_2 = \|\widehat{f}\|_2, \qquad f \in L^1(\mathbb{K}) \cap L^2(\mathbb{K}).
\]
Moreover, $\|\widehat{f}\|_\infty \leq \|f\|_1$ and the inversion formula holds:
\[
f(x,t) = \int_{\widehat{\mathbb{K}}} \widehat{f}(\lambda,m)\, \varphi_{(\lambda,m)}(x,t)\, d\gamma_\alpha(\lambda,m).
\]

Dilations play an important role in our analysis. For $r>0$, we define
\[
\delta_r(x,t) := (rx,r^2t)m, \qquad (\delta_r f)(x,t) := r^{-Q}\, f\!\left(\frac{x}{r},\frac{t}{r^2}\right), 
\]
where $Q=2\alpha+4$ is the homogeneous dimension of $\mathbb{K}$. It is easy to verify that
\[
\|\delta_r f\|_1 = \|f\|_1 \quad \text{and} \quad \widehat{\delta_r f}(\lambda,m) = \widehat{f}(r^2\lambda,m).
\]

\medskip

Throughout the paper, the notation $A \lesssim B$ means that there exists a constant $c>0$, depending only on the parameters $a,b,p$ and $\alpha$, such that $A \leq cB$.


\section{Proof of the Main Result}\label{sec3}

The main objective of this section is to establish our main result, Theorem \ref{1.5}. The proof naturally divides into three cases according to the range of $p$, namely $p=1$, $1 < p < 2$ and $p = 2$. We treat each case separately by formulating and proving a dedicated theorem. Let us begin with the case $p=1$.

\subsection{Case $p=1$}

In this case, we establish the desired inequality in the endpoint space $L^1(\mathbb{K})$. The strategy is as follows: we first observe that the inequality is invariant under dilation and scaling. This allows us to normalize the problem, thereby reducing the proof to showing a uniform lower bound on the Fourier side. Once this normalization is in place, we exploit the $L^2$ theory, together with the Plancherel identity and appropriate localization arguments, to obtain the necessary estimate. More precisely, we have the following theorem.

\begin{thm}\label{case1}
Let $a > 0$ and $b>Q/2$. Then, for every $f \in L^1(\mathbb{K})$, we have
\begin{equation}\label{3.1.1}
    \|f\|_1^{a+b} \leq C_{a,b,\alpha}\,\|\, |(x,t)|^a f\|_1^b \, \|\, |(\lambda,m)|^{b/2} \hat{f}\|_\infty^a,
\end{equation}
where $C_{a,b,\alpha}$ is a positive constant depending only on the parameters.
\end{thm}

\begin{proof}
We first note that \eqref{3.1.1} is invariant under dilation and scaling. Indeed, if 
\[
g(x,t) = c\,\delta_r(f)(x,t), \qquad c,r>0,
\]
then
\[
\|g\|_1^{a+b} = c^{a+b}\|f\|_1^{a+b}, \qquad 
\|\,|(x,t)|^a g\|_1^b = c^b r^{ab}\,\|\,|(x,t)|^a f\|_1^b,
\]
and
\[
\|\,|(\lambda,m)|^{b/2}\hat g\|_\infty^a = c^a r^{-ab}\,\|\,|(\lambda,m)|^{b/2}\hat f\|_\infty^a.
\]
The powers of $c$ and $r$ cancel and therefore the inequality holds for $f$ if and only if it holds for $g$. By a suitable choice of $c$ and $r$, we may normalize so that
\[
\|f\|_1 = 1, 
\qquad 
\|\, |(x,t)|^a f\|_1 = 1.
\]
Under this normalization, the inequality reduces to proving that
\begin{equation}\label{3.1.2}
    \|\, |(\lambda,m)|^{b/2}\hat f\|_\infty \;\geq\; C_{a,b,\alpha}.
\end{equation}

If $\|\,|(\lambda,m)|^{b/2}\hat f\|_\infty \geq 1$, the claim is trivial. We therefore assume
\[
\|\,|(\lambda,m)|^{b/2}\hat f\|_\infty \leq 1,
\]
so that
\begin{equation}\label{min}  
|\hat f(\lambda,m)| \leq \min\left\{1,\;\frac{1}{|(\lambda,m)|^{b/2}}\right\}.
\end{equation}

Let $E_r := \{(\lambda,m)\in \widehat{\mathbb{K}}: |(\lambda,m)|<r\}$. By \cite[Remarks~3.4(1)]{AR},
\begin{equation}\label{3.1.3}
\gamma_\alpha(E_r) = 2\frac{r^{\alpha+2}}{\alpha+2} \sum_{m \geq 0} \frac{L_m^\alpha(0)}{(4m+2\alpha+2)^{\alpha+2}} < \infty.
\end{equation}
Consequently, it follows from \eqref{min} and \eqref{3.1.3} that
\begin{align*}
    \int_{\widehat{\mathbb{K}}} |\hat f|^2\, d\gamma_\alpha 
&= \int_{E_1} |\hat f|^2\, d\gamma_\alpha + \int_{E_1^c} |\hat f|^2\, d\gamma_\alpha \\
&\leq \gamma_\alpha(E_1) +  \frac{2}{b-\alpha-2} \sum_{m \geq 0} \frac{L_m^\alpha(0)}{(4m+2\alpha+2)^{\alpha+2}} < \infty ,
\end{align*}
since $b>Q/2=\alpha+2$. Thus, $f\in L^2(\mathbb{K})$.

Consider $B_s:=\{(x,t)\in \mathbb{K}: |(x,t)|<s\}$. From the normalization,
\[
\int_{B_s} |f|\, dm_\alpha = \int_{\mathbb{K}} |f|\, dm_\alpha - \int_{B_s^c} |f|\, dm_\alpha \geq 1-s^{-a}.
\]
Hence for suitable $s_0>0$ with $0<1-s_0^{-a}<1$ and by applying the H\"{o}lder's inequality, we have
\[
\int_{B_{s_0}} |f|^2\, dm_\alpha \;\geq\; \frac{(1-s_0^{-a})^2}{s_0^Q \Omega_\alpha},
\]
where $\Omega_\alpha$ is the volume of the unit ball $B_1$. By the Plancherel identity, it follows that
\begin{equation}\label{3.1.4}
    \int_{\widehat{\mathbb{K}}} |\hat f|^2\, d\gamma_\alpha 
\;\geq\; \frac{(1-s_0^{-a})^2}{s_0^Q \Omega_\alpha}.
\end{equation}

On the other hand, for $r>0$,
\begin{equation}\label{3.1.5}
    \int_{E_r^c} |\hat f|^2\, d\gamma_\alpha \leq \frac{2}{b-\alpha-2}\,\frac{1}{r^{b-\alpha-2}} \sum_{m \geq 0}\frac{L_m^\alpha(0)}{(4m+2\alpha+2)^{\alpha+2}}.
\end{equation}
Since $|\hat{f}(\lambda,m)| \leq 1,$ it follows from \eqref{3.1.4} and \eqref{3.1.5} that there exists $r_1>0$ such that
\[
\int_{E_{r_1}} |\hat f|\, d\gamma_\alpha \;\geq\; \int_{E_{r_1}} |\hat{f}|^2 \, d\gamma_\alpha \;\geq\; \frac{1}{2}\,\frac{(1-s_0^{-a})^2}{s_0^Q \Omega_\alpha}.
\]
Further, using \eqref{3.1.3}, there exists $r_2>0$ with
\[
\int_{r_2 \leq |(\lambda,m)|<r_1} |\hat f|\, d\gamma_\alpha \;\geq\; \frac{1}{4}\,\frac{(1-s_0^{-a})^2}{s_0^Q \Omega_\alpha}.
\]

Finally, using the bound on $\hat f$, we obtain
\[
\int_{r_2 \leq |(\lambda,m)|<r_1} |\hat f|\, d\gamma_\alpha
\leq \frac{\gamma_\alpha(E_{r_1})}{r_2^{b/2}} \|\,|(\lambda,m)|^{b/2}\hat f\|_\infty.
\]
Combining the above estimates yields
\[
\|\,|(\lambda,m)|^{b/2}\hat f\|_\infty \;\geq\; \frac{r_2^{b/2}(1-s_0^{-a})^2}{4\,\gamma_\alpha(E_{r_1})\, s_0^Q \Omega_\alpha}.
\]
Thus \eqref{3.1.2} holds with
\[
C_{a,b,\alpha} := \frac{r_2^{b/2}(1-s_0^{-a})^2}{4\,\gamma_\alpha(E_{r_1})\, s_0^Q \Omega_\alpha} > 0,
\]
thereby completing the proof.
\end{proof}

\subsection{Case: $1 < p < 2$}  

The uncertainty inequality continues to hold in this range, but its proof requires a careful blend of the ideas used in the endpoint case $p=1$ and additional integrability arguments. In contrast to the case $p=1,$ where the conjugate exponent is $p' = \infty$, we now have $p'= p/(p-1)$ which is finite and greater than $2$. This allows us to exploit Hölder’s inequality on both the physical and Fourier sides in a more refined manner.  

A central step in our approach is to establish that under an appropriate integrability condition involving a weighted Fourier norm, the functions in $L^p(\mathbb{K})$ must in fact belong to $L^2(\mathbb{K})$. This phenomenon, which may be interpreted as a kind of ``regularization" property, is captured by the following lemma.  

\begin{lem}\label{mainlem}
Let $1 < p < 2$ and $2 < p'= p/(p-1) < \infty$. Suppose that $b > Q(1/p - 1/2)$. If $f \in L^p(\mathbb{K})$ and  
\[
\|\,|(\lambda,m)|^{b/2}\hat{f}\|_{p'} < \infty,
\]  
then $f \in L^2(\mathbb{K}).$
\end{lem}  

\begin{proof}  
The idea is to show that $\hat{f} \in L^2(\widehat{\mathbb{K}}),$ since by \cite[Theorem II.2]{NT}, this condition is equivalent to $f \in L^2(\mathbb{K}).$  

We begin by splitting the Fourier domain into two regions: one where $|(\lambda,m)|$ is small and another where it is large. For $r>0,$ define  
\[
E_r = \{(\lambda,m)\in \widehat{\mathbb{K}} : |(\lambda,m)| < r\}.
\]  
    Then  
    \[
    \int_{\widehat{\mathbb{K}}} |\hat{f}|^2 \, d\gamma_\alpha 
    = \underbrace{\int_{E_r} |\hat{f}|^2 \, d\gamma_\alpha}_{I_1} 
    + \underbrace{\int_{E_r^c} |\hat{f}|^2 \, d\gamma_\alpha}_{I_2}.
    \]  

    \medskip
 
    By H\"{o}lder’s inequality,  
    \begin{align*}
            I_1 &= \int_{E_r} |\hat{f}(\lambda,m)|^2 \, d\gamma_\alpha(\lambda,m) \\
            &\leq \left( \int_{E_r} |\hat{f}(\lambda,m)|^{p'} \, d\gamma_\alpha(\lambda,m) \right)^{\frac{2}{p'}} \, \left( \int_{E_r} d\gamma_\alpha(\lambda,m) \right)^{1 - \frac{2}{p'}} \\
            &\leq \|\hat{f}\|_{p'}^2 \, (\gamma_\alpha(E_r))^{1 - \frac{2}{p'}}.
        \end{align*}
    Since $1<p < 2,$ the Hausdorff–Young inequality \cite[Theorem 3.4]{VM} gives  
    $$\|\hat{f}\|_{p'} \lesssim \|f\|_p,$$
        for $f \in L^p(\mathbb{K})$. Thus,
        $$I_1 \lesssim \|f\|_p^2 \, (r^{\alpha+2})^{1 - \frac{2}{p'}} = \|f\|_p^2 \, r^{\frac{Q}{2}\frac{2-p}{p}}< \infty.$$

    \medskip
 
    For the complementary region $E_r^c,$ we consider the integral $I_2$ and again, by applying the H\"{o}lder's inequality, we have
        \begin{align*}
            I_2 &= \int_{E_r^c} |\hat{f}(\lambda,m)|^2 \, d\gamma_\alpha(\lambda,m) \\
            &= \int_{E_r^c} \left( |(\lambda,m)|^{b/2} \, |\hat{f}(\lambda,m)| \right)^2 \, |(\lambda,m)|^{-b} \, d\gamma_\alpha(\lambda,m) \\
            &\leq \left( \int_{E_r^c} \left( |(\lambda,m)|^{b/2} \, |\hat{f}(\lambda,m)| \right)^{p'} d\gamma_\alpha(\lambda,m) \right)^\frac{2}{p'} \, \left( \int_{E_r^c} |(\lambda,m)|^\frac{-b p'}{p' - 2} \, d\gamma_\alpha(\lambda,m) \right)^{1 - \frac{2}{p'}} \\
            &\leq \||(\lambda,m)|^{b/2} \, \hat{f}\|_{p'}^2 \,\, (I_3)^{1 - \frac{2}{p'}},
        \end{align*}
        where 
        $$I_3 := \int_{E_r^c} |(\lambda,m)|^\frac{-b p'}{p' - 2} \, d\gamma_\alpha(\lambda,m).$$
       Since $b > Q(1/p - 1/2),$ that is, $\alpha + 2 - \frac{b p'}{p' - 2} < 0,$ we obtain
        \begin{align*}
            I_3  &= 2 \sum_{m \geq 0} \frac{L_m^\alpha(0)}{(4m+2\alpha+2)^{\frac{b p'}{p' - 2}}} \int_{\frac{r}{4m+2\alpha+2}}^\infty \lambda^{\frac{-b p'}{p' - 2}} \,  \, \lambda^{\alpha+1} \, d\lambda \\
            &= \frac{2}{\frac{b p'}{p' - 2} - \alpha - 2} \, \frac{1}{r^{\frac{b p'}{p' - 2}-\alpha - 2}} \, \sum_{m \geq 0} \frac{L_m^\alpha(0)}{(4m+2\alpha+2)^{\alpha + 2}} < \infty.
        \end{align*}

    \medskip
    \noindent Combining the two estimates, both $I_1$ and $I_2$ are finite, which implies $\hat{f} \in L^2(\widehat{\mathbb{K}}).$ Therefore, $f \in L^2(\mathbb{K}),$ which completes the proof.  
    \end{proof}  

    \medskip

    With Lemma \ref{mainlem} in place, we can now establish the main result for the case $1<p<2.$  

    \begin{thm}\label{case2}
    Let $1<p<2$ and $2<p'= p/(p-1) <\infty.$ Suppose that $a > 0$ and $b > Q(1/p - 1/2).$ Then, for every $f \in L^p(\mathbb{K}),$ we have
    \begin{equation}\label{3.2.1}
    \|f\|_p^{a+b} \leq C_{a,b,p,\alpha} \,\|\,|(x,t)|^a f\|_p^b \,\|\,|(\lambda,m)|^{b/2}\hat{f}\|_{p'}^a,
    \end{equation}  
    where $C_{a,b,p,\alpha}$ is a positive constant depending only on the parameters.
    \end{thm}  

    \begin{proof}
    As in the case $p=1$, we first prove that the desired inequality is invariant 
    with respect to both dilation and scaling by a constant. 
    For $c,r > 0$ and $f \in L^p(\mathbb{K})$, let $g := c \, \delta_r(f)$. 
    With simple computations, it is easy to verify the following

     $$\|g\|_p^{a+b} = c^{a+b} \, r^{\frac{-Q(a+b)}{p'}} \, \|f\|_p^{a+b},$$
        $$\| \, |(x,t)|^a \, g \|_p^b = c^b \, r^{ab - \frac{Q b}{p'}} \, \| \, |(x,t)|^a \, f \|_p^b$$
        and
        $$\| \, |(\lambda,m)|^{b/2} \, \hat{g}\|_{p'}^a = c^a \, r^{-ab - \frac{Qa}{p'}} \, \| \, |(\lambda,m)|^{b/2} \, \hat{f}\|_{p'}^a.$$
        This proves the invariance of the inequality, so without loss of generality, we may assume that 
        $$\|f\|_p = 1 = \| \, |(x,t)|^a \, f \|_p.$$
        Therefore, in order to establish \eqref{3.2.1}, it suffices to prove that there exists a positive constant $C_{a,b,p,\alpha}$ such that
        \begin{equation}\label{3.2.2}
            \| \, |(\lambda,m)|^{b/2} \, \hat{f}\|_{p'} \geq C_{a,b,p,\alpha}.
        \end{equation}
        Following the same arguments as before, without loss of generality, we may assume that
        $$\| \, |(\lambda,m)|^{b/2} \, \hat{f}\|_{p'} \leq 1.$$
        Since $b > Q(1/p - 1/2)$, by Lemma \ref{mainlem} and the Plancherel identity, it follows that $f \in L^2(\mathbb{K})$ and
        $$\|f\|_2 = \|\hat{f}\|_2 \lesssim r^{\frac{Q}{2}\frac{2-p}{p}} + \left( r^{\frac{Q}{2}-\frac{b p}{2 - p}} \right)^{\frac{2 - p}{p}} = r^{\frac{Q}{2} \left( \frac{2 - p}{p} \right)} (1 + r^{-b}),$$
        for any $r > 0.$
        
        Now, for $s > 0,$ let $B_s := \{(x,t) \in \mathbb{K}: |(x,t)| < s\}.$ Given our initial assumption, we have
        $$1 = \int_\mathbb{K} |(x,t)|^{ap} \, |f(x,t)|^p \, dm_\alpha(x,t) \geq s^{ap} \int_{B_s^c} |f(x,t)|^p \, dm_\alpha(x,t).$$
        This implies that
        $$\int_{B_s} |f(x,t)|^p \, dm_\alpha(x,t) = \int_\mathbb{K} |f(x,t)|^p \, dm_\alpha(x,t) - \int_{B_s^c} |f(x,t)|^p \, dm_\alpha(x,t) \geq 1 - s^{-{ap}}.$$
        Choosing $s = s_0 > 0$ such that $0 < 1 - {s_0}^{-ap} < 1$ and by applying the H\"{o}lder's inequality, we obtain
        \begin{align*}
            1-{s_0}^{-ap} &\leq \int_{B_{s_0}} |f(x,t)|^p \, dm_\alpha(x,t) \\
            &\leq \left( \int_{B_{s_0}} |f(x,t)|^2 dm_\alpha(x,t) \right)^\frac{p}{2} \, \left( \int_{B_{s_0}} dm_\alpha(x,t) \right)^{1-\frac{p}{2}} \\
            &\lesssim \left( \int_{B_{s_0}} |f(x,t)|^2 dm_\alpha(x,t) \right)^\frac{p}{2} \, {s_0}^{Q \left(1 - \frac{p}{2} \right)}.
        \end{align*}
        Now, by the Plancherel identity,
        \begin{align*}
            \int_{\widehat{\mathbb{K}}} |\hat{f}(\lambda,m)|^2 \, d\gamma_\alpha(\lambda,m) &= \int_\mathbb{K} |f(x,t)|^2 \, dm_\alpha(x,t) \\
            &\geq \int_{B_{s_0}} |f(x,t)|^2 \, dm_\alpha(x,t) \\
            &\gtrsim (1 - {s_0}^{-{ap}})^{\frac{2}{p}}{s_0}^{Q \left(1 - \frac{2}{p} \right)}.
        \end{align*}
        Further, for $r > 0$ and $b > Q(1/p - 1/2),$ observe that
        \begin{align*}
            \int_{E_r} |\hat{f}(\lambda,m)|^2 \, d\gamma_\alpha(\lambda,m) &= \int_{\widehat{\mathbb{K}}} |\hat{f}(\lambda,m)|^2 \, d\gamma_\alpha(\lambda,m) - \int_{E_r^c} |\hat{f}(\lambda,m)|^2 \, d\gamma_\alpha(\lambda,m) \\
            &\gtrsim (1 - {s_0}^{-{ap}})^{\frac{2}{p}} \, {s_0}^{Q \left(1 - \frac{2}{p} \right)} - \int_{E_r^c} |\hat{f}(\lambda,m)|^2 \, d\gamma_\alpha(\lambda,m).
        \end{align*}
        From the proof of Lemma \ref{mainlem}, we have
        \begin{align*}
            \int_{E_r} |\hat{f}(\lambda,m)|^2 \, d\gamma_\alpha(\lambda,m) &\gtrsim (1 - {s_0}^{-{ap}})^{\frac{2}{p}} \, {s_0}^{Q \left(1 - \frac{2}{p} \right)} - \left( r^{\frac{Q}{2}-\frac{b p}{2 - p}} \right)^{\frac{2 - p}{p}}.
        \end{align*}
        Choose $r = r_1 > 0$ (depending on $a,b,p$ and $\alpha$) large enough such that 
        \begin{align*}
            \int_{E_r} |\hat{f}(\lambda,m)|^2 \, d\gamma_\alpha(\lambda,m) &\gtrsim \frac{1}{2}(1 - {s_0}^{-{ap}})^{\frac{2}{p}}{s_0}^{Q \left(1 - \frac{2}{p} \right)}.
        \end{align*}
        Then, for any $r > 0,$ it follows that
        \begin{align*}
            \int_{r \leq |(\lambda,m)| < r_1} |\hat{f}(\lambda,m)|^2 \, d\gamma_\alpha(\lambda,m) &= \int_{E_{r_1}} |\hat{f}(\lambda,m)|^2 \, d\gamma_\alpha(\lambda,m) - \int_{E_{r}} |\hat{f}(\lambda,m)|^2 \, d\gamma_\alpha(\lambda,m) \\
            &\gtrsim \frac{1}{2}(1 - {s_0}^{-{ap}})^{\frac{2}{p}} \, {s_0}^{Q \left(1 - \frac{2}{p} \right)} - \int_{E_r} |\hat{f}(\lambda,m)|^2 \, d\gamma_\alpha(\lambda,m).
        \end{align*}
        Again, from the proof of Lemma \ref{mainlem}, we obtain
        \begin{align*}
            \int_{r \leq |(\lambda,m)| < r_1} |\hat{f}(\lambda,m)|^2 \, d\gamma_\alpha(\lambda,m) \gtrsim \frac{1}{2}(1 - {s_0}^{-{ap}})^{\frac{2}{p}} \, {s_0}^{Q \left(1 - \frac{2}{p} \right)} - r^{\frac{Q}{2}\frac{2-p}{p}}.
        \end{align*}
        Choose $r = r_2 > 0$ (depending on $a,b,p$ and $\alpha$) small enough such that
        \begin{align*}
            I := \int_{r_2 \leq |(\lambda,m)| < r_1} |\hat{f}(\lambda,m)|^2 \, d\gamma_\alpha(\lambda,m) \gtrsim \frac{1}{4}(1 - {s_0}^{-{ap}})^{\frac{2}{p}} \, {s_0}^{Q \left(1 - \frac{2}{p} \right)}.
        \end{align*}
        Now, by applying the H\"{o}lder's inequality, we have
        \begin{align*}
            I &\leq \left( \int_{r_2 \leq |(\lambda,m)| < r_1} |\hat{f}(\lambda,m)|^{p'} \, d\gamma_\alpha(\lambda,m) \right)^{\frac{2}{p'}} \, \left( \int_{r_2 \leq |(\lambda,m)| < r_1} d\gamma_\alpha(\lambda,m) \right)^{1 - \frac{2}{p'}} \\
            &\leq \left( \int_{r_2 \leq |(\lambda,m)| < r_1} |\hat{f}(\lambda,m)|^{p'} \, d\gamma_\alpha(\lambda,m) \right)^{\frac{2}{p'}} \, \left( \gamma_\alpha(E_{r_1}) \right)^{1 - \frac{2}{p'}}.
        \end{align*}
        Therefore,
        \begin{align*}
            \int_{r_2 \leq |(\lambda,m)| < r_1} |\hat{f}(\lambda,m)|^{p'} \, d\gamma_\alpha(\lambda,m) &\geq (I)^{\frac{p'}{2}} \, \left( \gamma_\alpha(E_{r_1}) \right)^{\frac{2 - p'}{2}} \\
            &\gtrsim \left( \frac{1}{4}(1 - {s_0}^{-{ap}})^{\frac{2}{p}} \, {s_0}^{Q \left(1 - \frac{2}{p} \right)} \right)^{\frac{p'}{2}} \, \left( \gamma_\alpha(E_{r_1}) \right)^{\frac{2 - p'}{2}}.
        \end{align*}
        Finally, we have
        \begin{align*}
            \| \, |(\lambda,m)|^{b/2} \, \hat{f} \|_{p'}^{p'} &= \int_{\widehat{\mathbb{K}}} |(\lambda,m)|^{bp'/2} \, |\hat{f}(\lambda,m)|^{p'} \, d\gamma_\alpha(\lambda,m) \\
            &\geq \int_{r_2 \leq |(\lambda,m)| < r_1} |(\lambda,m)|^{bp'/2} \, |\hat{f}(\lambda,m)|^{p'} \, d\gamma_\alpha(\lambda,m) \\
            &\geq (r_2)^{bp'/2} \, \int_{r_2 \leq |(\lambda,m)| < r_1} |\hat{f}(\lambda,m)|^{p'} \, d\gamma_\alpha(\lambda,m) \\
            &\gtrsim (r_2)^{bp'/2} \, \left( \frac{1}{4}(1 - {s_0}^{-{ap}})^{\frac{2}{p}} \, {s_0}^{Q \left(1 - \frac{2}{p} \right)} \right)^{\frac{p'}{2}} \, \left( \gamma_\alpha(E_{r_1}) \right)^{\frac{2 - p'}{2}}.
        \end{align*}
        This implies that
        \begin{align*}
            \| \, |(\lambda,m)|^{b/2} \, \hat{f} \|_{p'} \gtrsim \left( (r_2)^{bp'/2} \, \left( \frac{1}{4}(1 - {s_0}^{-{ap}})^{\frac{2}{p}} \, {s_0}^{Q \left(1 - \frac{2}{p} \right)} \right)^{\frac{p'}{2}} \, \left( \gamma_\alpha(E_{r_1}) \right)^{\frac{2 - p'}{2}} \right)^{\frac{1}{p'}}.
        \end{align*}
       Thus, there exists a positive constant $C_{a,b,p,\alpha}$ such that
    \[
    \|\, |(\lambda,m)|^{b/2} \, \hat{f}\,\|_{p'} \geq C_{a,b,p,\alpha},
    \]
    thereby establishing \eqref{3.2.2}. 
    This completes the proof of the theorem.
    \end{proof}

    The following corollary establishes an HPW-type uncertainty inequality for the Schwartz class functions on $\mathbb{K}$ in the range $2 < p < \infty$. It is a direct consequence of Theorem \ref{case2}, applied with the conjugate exponent $p'$, together with the Hausdorff--Young inequality.

    \begin{cor}
        Let $2 < p < \infty$ and $1 < p'= p/(p-1) < 2$. Suppose that $a > 0$ and $b > Q(1/2 - 1/p).$ Then, for every $f \in \mathcal{S}(\mathbb{K}),$ we have
    \begin{equation*}
    \|\hat{f}\|_p^{a+b} \leq C_{a,b,p,\alpha} \,\|\,|(x,t)|^a f\|_{p'}^b \,\|\,|(\lambda,m)|^{b/2}\hat{f}\|_{p}^a,
    \end{equation*}  
    where $C_{a,b,p,\alpha}$ is a positive constant depending only on the parameters.
    \end{cor}

    \subsection{Case: $p = 2$}

   We establish the classical Heisenberg--Pauli--Weyl uncertainty inequality on the Laguerre hypergroup. Our proof is entirely independent of heat kernel techniques and thus provides an alternative approach to the problem. Moreover, it strengthens the previously known result \cite[Theorem 4.1]{AR} by removing the assumption $a,b \geq 1$, thereby yielding the inequality in full generality. We begin by establishing an auxiliary lemma on exponential weight functions, which will serve as a fundamental tool in the proof of the main result.

    For each $s > 0,$ define the function $g_s: \widehat{\mathbb{K}} \to \mathbb{R}$ by
    $$g_s(\lambda,m) := e^{-s \, |(\lambda,m)|^2}.$$
    Clearly, $g_s \in L^\infty(\widehat{\mathbb{K}})$ for every $s > 0.$

    \begin{lem}\label{lem}
    Let $1 \leq p < \infty.$ Then, for every $s > 0,$ we have $g_s \in L^p(\widehat{\mathbb{K}})$ and $$\|g_s\|_{p} \lesssim s^{\frac{-Q}{4p}}.$$
    \end{lem}

    \begin{proof}
        For $s > 0,$ we have
        \begin{align*}
            \|g_s\|_p^p &= \int_{\widehat{\mathbb{K}}} |g_s(\lambda,m)|^p \, d\gamma_\alpha(\lambda,m)\\ &= \int_{\widehat{\mathbb{K}}} \left( e^{-s \, |(\lambda,m)|^2} \right)^p \, d\gamma_\alpha(\lambda,m) \\
            &= 2 \sum_{m \geq 0} L_m^\alpha(0) \int_0^\infty e^{-sp \, \lambda^2 (4m+2\alpha+2)^2} \, \lambda^{\alpha+1} \, d\lambda \\
            &= 2 \sum_{m \geq 0} L_m^\alpha(0) \left( \frac{1}{2} \frac{1}{s^{Q/4} \, p^{Q/4}} \, \Gamma\left(\frac{Q}{4}\right) \frac{1}{(4m+2\alpha+2)^{\alpha+2}} \right) \\
            &= \frac{s^{-Q/4}}{p^{Q/4}} \, \Gamma\left(\frac{Q}{4}\right) \sum_{m \geq 0} \frac{L_m^\alpha(0)}{(4m+2\alpha+2)^{\alpha+2}} < \infty.
        \end{align*}
        This implies that 
        $$g_s \in L^p(\widehat{\mathbb{K}}) \,\, \text{and} \,\, \|g_s\|_p \lesssim s^{\frac{-Q}{4p}}.$$
    \end{proof}

    We are now in a position to establish the classical Heisenberg--Pauli--Weyl uncertainty inequality on the Laguerre hypergroup.

    \begin{thm}\label{case3}
    Let $a,b  > 0$. Then, for every $f \in L^2(\mathbb{K})$, we have
    \begin{equation}\label{3.3.1}
        \| f \|_2^{\,a+b} \lesssim \| \, |(x,t)|^a f\|_2^b \, \| \, |(\lambda,m)|^{b/2} \hat{f}\|_2^a.
    \end{equation}
    \end{thm}

    \begin{proof}
        Let $s > 0$ and $f \in L^2(\mathbb{K})$. To prove \eqref{3.3.1}, we first establish the estimate
        \begin{equation}\label{3.3.2}
            \|g_s \, \hat{f}\|_2 \lesssim s^{-a/4} \, \| \, |(x,t)|^a \, f\|_{2}.
        \end{equation}
        For $r > 0,$ define $$B_r := \{(x,t) \in \mathbb{K}: |(x,t)| < r\}.$$ Observe that
        $$|f(x,t) \, \chi_{B_r^c}(x,t)| \leq r^{-a} \, |(x,t)|^a \, |f(x,t)|.$$
        Using this estimate together with the Plancherel identity, we obtain
        $$\|g_s \,  \widehat{f \, \chi_{B_r^c}} \|_2 \leq \|g_s\|_\infty \, \| \widehat{f \, \chi_{B_r^c}} \|_2 \lesssim \| f \, \chi_{B_r^c} \|_{2} \leq r^{-a} \, \| \, |(x,t)|^a \, f\|_{2}.$$
        On the other hand, using the estimate $\| \hat{f} \|_\infty \leq \|f\|_1$ together with Lemma \ref{lem} (with $p = 2$), we have
        $$\|g_s \,  \widehat{f \, \chi_{B_r}} \|_2 \leq \|g_s\|_2 \, \| \widehat{f \, \chi_{B_r}} \|_\infty \leq \|g_s\|_2 \, \|f \, \chi_{B_r}\|_1 \lesssim s^{-Q/8} \, \| \, |(x,t)|^{-a} \, \chi_{B_r}\|_2 \, \| \, |(x,t)|^a \, f\|_{2}.$$
        If $0 < a < Q/2$, then
        \begin{align*}
            \| \, |(x,t)|^{-a} \, \chi_{B_r}\|_2 &= \left( \int_{B_r} |(x,t)|^{-2a} \, dm_\alpha(x,t) \right)^{1/2} \lesssim r^{\frac{Q - 2a}{2}}. 
        \end{align*}
        Consequently,
        $$\|g_s \,  \widehat{f \, \chi_{B_r}} \|_2 \lesssim s^{-Q/8} \, r^{\frac{Q - 2a}{2}} \, \| \, |(x,t)|^a \, f\|_{2}.$$
        Combining the above estimates yields
        \begin{align*}
            \|g_s \, \hat{f}\|_2 &\leq \|g_s \, \widehat{f \, \chi_{B_r}}\|_2 + \|g_s \, \widehat{f \, \chi_{B_r^c}}\|_2 \\
            &\lesssim \left( s^{-Q/8} \, r^{\frac{Q - 2a}{2}} + r^{-a} \right) \| \, |(x,t)|^a \, f\|_{2}.
        \end{align*}
        Finally, choosing $r = s^{1/4} > 0,$ we obtain
        $$\|g_s \, \hat{f}\|_2 \lesssim s^{-a/4} \, \| \, |(x,t)|^a \, f\|_{2}.$$ 
        This establishes the desired estimate.

        To prove \eqref{3.3.1}, we may assume without loss of generality that $$\| \, |(x,t)|^a \, f\|_{2} \neq 0, \qquad \| \, |(\lambda,m)|^{b/2} \, \hat{f}\|_2 \neq 0,$$ and that both quantities are finite.

        \noindent By \eqref{3.3.2}, we have
        \begin{align*}
            \| \hat{f} \|_2 &\leq \|g_s \, \hat{f}\|_2 + \| (1-g_s) \, \hat{f}\|_2 \\
            &\lesssim s^{-a/4} \, \| \, |(x,t)|^a \, f\|_{2} + \|(1 - e^{-s \, |(\lambda,m)|^2}) \, \hat{f}\|_2 \\
            &= s^{-a/4} \, \| \, |(x,t)|^a \, f\|_{2} + s^{b/4} \, \| (s \, |(\lambda,m)|^2)^{-b/4} \, (1 - e^{-s \, |(\lambda,m)|^2}) \, |(\lambda,m)|^{b/2} \, \hat{f}\|_2.
        \end{align*}
        If $0 < b < 4$, then $$z^{-b/4} \, (1 - e^{-z}) \leq 1, \qquad z > 0,$$
        and therefore
        $$\| \hat{f} \|_2 \lesssim s^{-a/4} \, \| \, |(x,t)|^a \, f\|_{2} + s^{b/4} \, \| \, |(\lambda,m)|^{b/2} \, \hat{f}\|_2.$$
        Choosing $s > 0$ such that $$s^{1/4} = \left( \frac{a \, \| \, |(x,t)|^a \, f\|_{2}}{b \, \| \, |(\lambda,m)|^{b/2} \, \hat{f}\|_2} \right)^{\frac{1}{a+b}} > 0,$$
        yields
        $$\| \hat{f} \|_2 \lesssim \| \, |(x,t)|^a \, f\|_{2}^{\frac{b}{a+b}} \, \| \, |(\lambda,m)|^{b/2} \, \hat{f}\|_2^{\frac{a}{a+b}}.$$
        Invoking the Plancherel identity, we obtain
        $$\| f \|_2^{a+b} \lesssim \| \, |(x,t)|^a \, f\|_{2}^b \, \| \, |(\lambda,m)|^{b/2} \, \hat{f}\|_2^a,$$
        whenever $0 < a < Q/2$ and $0 < b < 4$.

        It remains to treat the case $a \geq Q/2$. Let $a' < Q/2$. Then, for any $\epsilon > 0$, 
        $$\left(\frac{|(x,t)|}{\sqrt{\epsilon}}\right)^{a'} \leq 1 + \left(\frac{|(x,t)|}{\sqrt{\epsilon}}\right)^{a},$$
        which gives
        $$\| \, |(x,t)|^{a'} \, f\|_{2} \leq \epsilon^{a'/2} \, \|f\|_2 + \epsilon^{(a'-a)/2} \, \| \, |(x,t)|^a \, f\|_{2}.$$
        Choosing 
        $$\epsilon = \left( \frac{(a-a') \, \| \, |(x,t)|^a \, f\|_{2}}{a' \, \|f\|_2} \right)^{2/a},$$
        which minimizes the right-hand side, we obtain
        $$\| \, |(x,t)|^{a'} \, f\|_{2} \lesssim \|f\|_2^{(a-a')/a} \, \| \, |(x,t)|^a \, f\|_{2}^{a'/a}.$$
        Combining this estimate with \eqref{3.3.1}, applied with the exponent $a'$, we conclude that
        $$\| f \|_2^{a+b} \lesssim \| \, |(x,t)|^a \, f\|_{2}^b \, \| \, |(\lambda,m)|^{b/2} \, \hat{f}\|_2^a,$$
        for all $a > 0$ and $0 < b < 4$.

        Finally, the case $b \geq 4$ follows by an analogous H\"{o}lder interpolation argument, which shows that the above inequality remains valid for all $b \geq 4$. Consequently, \eqref{3.3.1} holds for all $a,b > 0$, thereby completing the proof.
    \end{proof}
    
    \section*{Concluding Remarks}
The main contributions of this paper are summarized as follows.

\begin{itemize}
\item We established the first $L^p$-Heisenberg--Pauli--Weyl uncertainty inequalities on the Laguerre hypergroup for the full range $1\le p\le2$, thereby extending Xiao's Euclidean $L^p$ theory to the hypergroup setting.

\item We proved a refined $L^2$-Heisenberg--Pauli--Weyl uncertainty inequality that removes the restrictions $a,b\ge1$ imposed in the earlier work of Atef. In particular, we showed that the inequality remains valid for all positive exponents $a,b>0$, demonstrating that the previous restrictions arise from the heat kernel approach rather than from the inequality itself.

\item Unlike the Euclidean setting, the Fourier analysis on the Laguerre hypergroup is carried out with respect to the mixed discrete--continuous spectral variables $(\lambda,m)$ and the associated Plancherel measure $d\gamma_\alpha(\lambda,m)$. Consequently, the proof require estimates for spectral balls, weighted Fourier norms and scaling properties adapted to the Fourier--Laguerre transform, thereby isolating the analytical features specific to the hypergroup framework.

\item Our proofs are entirely Fourier-analytic and rely on the Fourier--Laguerre transform, dilation and scaling invariance, the Hausdorff--Young inequality and the Plancherel identity, completely avoiding heat kernel techniques.

\item The present work provides a unified framework for Heisenberg--Pauli--Weyl uncertainty inequalities on the Laguerre hypergroup and highlights the connections between Euclidean harmonic analysis, the Heisenberg group and commutative hypergroups.
\end{itemize}

 \section*{Acknowledgement}
    The first author gratefully acknowledges the Indian Institute of Technology Delhi for providing the Institute Assistantship.

    \section*{Data Availability} 
    Data sharing does not apply to this article as no datasets were generated or analyzed during the current study.

    \section*{Competing Interests}
    The authors declare that they have no competing interests. 
    
	\bibliographystyle{acm}
	\bibliography{ref}

\end{document}